\documentclass[11pt,twoside]{article}
\usepackage{amsmath, amsthm, amscd, amsfonts, amssymb, graphicx, color}

\date{}

\begin{document}

\centerline{}

\centerline {\Large{\bf A few fixed point theorems in linear $n$-normed space}}

\newcommand{\mvec}[1]{\mbox{\bfseries\itshape #1}}
\centerline{}
\centerline{\textbf{Prasenjit Ghosh}}
\centerline{Department of Pure Mathematics, University of Calcutta}
\centerline{35, Ballygunge Circular Road, Kolkata, 700019, West Bengal, India}
\centerline{e-mail: prasenjitpuremath@gmail.com}

\centerline{}
\centerline{\textbf{T. K. Samanta}}
\centerline{Department of Mathematics, Uluberia College}
\centerline{Uluberia, Howrah, 711315,  West Bengal, India}
\centerline{e-mail: mumpu$_{-}$tapas5@yahoo.co.in}

\newtheorem{Theorem}{\quad Theorem}[section]

\newtheorem{definition}[Theorem]{\quad Definition}

\newtheorem{theorem}[Theorem]{\quad Theorem}

\newtheorem{remark}[Theorem]{\quad Remark}

\newtheorem{corollary}[Theorem]{\quad Corollary}

\newtheorem{note}[Theorem]{\quad Note}

\newtheorem{lemma}[Theorem]{\quad Lemma}

\newtheorem{example}[Theorem]{\quad Example}

\newtheorem{result}[Theorem]{\quad Result}
\newtheorem{conclusion}[Theorem]{\quad Conclusion}

\newtheorem{proposition}[Theorem]{\quad Proposition}

\begin{abstract}
\textbf{\emph{Banach's fixed point theorem in linear\;$n$-normed space is being developed.\,Also, we present several theorems on fixed points in linear\;$n$-normed space.}}
\end{abstract}


\section{Introduction}

\smallskip\hspace{.6 cm}
The Banach fixed point theorem concerns certain mappings of a complete metric space into itself.\,It states sufficient conditions for the existence and uniqueness of a fixed point.\,Also, the theorem gives an iterative process by which we can obtain approximations to the fixed point and error bounds.\,This theorem has important applications to finding the unique solution of linear algebraic equations, differential equations, integral equations and as well as to implicit function theorem.

The idea of linear 2-normed space was first introduced by S. Gahler \cite{Gahler} and thereafter the geometric structure of linear 2-normed spaces was developed by A. White, Y. J. Cho, R. W. Freese \cite{White,Freese}.\;In recent times, some important results in classical normed spaces have been proved into 2-norm setting by many researchers.\,The concept of \,$2$-Banach space is briefly discussed in \cite{White}.\;H.\,Gunawan and Mashadi \cite{Mashadi} developed the generalization of a linear $2$-normed space for \,$n \,\geq\, 2$.\,Some fundamental results of classical normed space with respect to \,$b$-linear functional in linear\;$n$-normed space  have been studied by P. Ghosh and T. K. Samanta \cite{Prasenjit}.\,Also they have studied the reflexivity of linear \,$n$-normed space with respect to \,$b$-linear functional in \cite{K} and slow convergence of sequence of \,$b$-linear functionals in linear \,$n$-normed space in \cite{KK}.

In this paper, we first define the contraction mapping in linear \,$n$-normed space and then present the Banach fixed point theorem in linear \,$n$-normed space.\,Next, we discuss some further theorems on fixed points of operators in linear \,$n$-normed space.

\section{Preliminaries}
\smallskip\hspace{.6 cm}
In this section, we give some necessary definitions.

\begin{definition}\cite{Mashadi}
Let \,$X$\, be a linear space over the field \,$ \mathbb{K}$, where \,$ \mathbb{K} $\, is the real or complex numbers field with \,$\text{dim}\,X \,\geq\, n$, where \,$n$\, is a positive integer.\;A real valued function \,$\left \|\,\cdot \,,\, \cdots \,,\, \cdot \,\right \| \,:\, X^{\,n} \,\to\, \mathbb{R}$\, is called an n-norm on \,$X$\, if
\begin{itemize}
\item[(N1)]$\left\|\,x_{\,1} \,,\, x_{\,2} \,,\, \cdots \,,\, x_{\,n}\,\right\| \,=\,0$\, if and only if \,$x_{\,1},\, \cdots,\, x_{\,n}$\, are linearly dependent,
\item[(N2)]$\left\|\,x_{\,1} \,,\, x_{\,2} \,,\, \cdots \,,\, x_{\,n}\,\right\|$\; is invariant under permutations of \,$x_{\,1},\, x_{\,2},\, \cdots,\, x_{\,n}$,
\item[(N3)]$\left\|\,\alpha\,x_{\,1} \,,\, x_{\,2} \,,\, \cdots \,,\, x_{\,n}\,\right\| \,=\, |\,\alpha\,|\, \left\|\,x_{\,1} \,,\, x_{\,2} \,,\, \cdots \,,\, x_{\,n}\,\right\|\; \;\;\forall \;\; \alpha \,\in\, \mathbb{K}$,
\item[(N4)]$\left\|\,x \,+\, y \,,\, x_{\,2} \,,\, \cdots \,,\, x_{\,n}\,\right\| \,\leq\, \left\|\,x \,,\, x_{\,2} \,,\, \cdots \,,\, x_{\,n}\,\right\| \,+\,  \left\|\,y \,,\, x_{\,2} \,,\, \cdots \,,\, x_{\,n}\,\right\|$
\end{itemize}
hold for all \,$x,\, y,\, x_{\,1},\, x_{\,2},\, \cdots,\, x_{\,n} \,\in\, X$.\;The pair \,$\left(\,X,\, \left \|\,\cdot,\, \cdots,\, \cdot \,\right \| \,\right)$\; is then called a linear n-normed space. 
\end{definition}

Throughout this paper, \,$X$\, will denote linear\;$n$-normed space over the field \,$\mathbb{K}$\, of complex or real numbers, associated with the $n$-norm \,$\|\,\cdot \,,\, \cdots \,,\, \cdot\,\|$.

\begin{definition}\cite{Mashadi}
A sequence \,$\{\,x_{\,k}\,\} \,\subseteq\, X$\, is said to converge to \,$x \,\in\, X$\; if 
\[\lim\limits_{k \to \infty}\,\left\|\,x_{\,k} \,-\, x \,,\, e_{\,2} \,,\, \cdots \,,\, e_{\,n} \,\right\| \,=\, 0\]
for every \,$ e_{\,2},\, \cdots,\, e_{\,n} \,\in\, X$\, and it is called a Cauchy sequence if 
\[\lim\limits_{l \,,\, k \to \infty}\,\left\|\,x_{\,l} \,-\, x_{\,k} \,,\, e_{\,2} \,,\, \cdots \,,\, e_{\,n}\,\right\| \,=\, 0\]
for every \,$ e_{\,2},\, \cdots,\, e_{\,n} \,\in\, X$.\;The space \,$X$\, is said to be complete or n-Banach space if every Cauchy sequence in this space is convergent in \,$X$.
\end{definition}

\begin{definition}\cite{Soenjaya}
We define the following open and closed ball in \,$X$: 
\[B_{\,\{\,e_{\,2} \,,\, \cdots \,,\, e_{\,n}\,\}}\,(\,a \,,\, \delta\,) \,=\, \left\{\,x \,\in\, X \,:\, \left\|\,x \,-\, a \,,\, e_{\,2} \,,\, \cdots \,,\, e_{\,n}\,\right\| \,<\, \delta \,\right\}\;\text{and}\]
\[B_{\,\{\,e_{\,2} \,,\, \cdots \,,\, e_{\,n}\,\}}\,[\,a \,,\, \delta\,] \,=\, \left\{\,x \,\in\, X \,:\, \left\|\,x \,-\, a \,,\, e_{\,2} \,,\, \cdots \,,\, e_{\,n}\,\right\| \,\leq\, \delta\,\right\},\hspace{.5cm}\]
where \,$a,\, e_{\,2},\, \cdots,\, e_{\,n} \,\in\, X$\, and \,$\delta$\, be a positive number.
\end{definition}

\begin{definition}\cite{Soenjaya}
A subset \,$G$\, of \,$X$\, is said to be open in \,$X$\, if for all \,$a \,\in\, G $, there exist \,$e_{\,2},\, \cdots,\, e_{\,n} \,\in\, X $\, and \, $\delta \,>\, 0 $\; such that \,$B_{\,\{\,e_{\,2} \,,\, \cdots \,,\, e_{\,n}\,\}}\,(\,a \,,\, \delta\,) \,\subseteq\, G$.
\end{definition}

\begin{definition}\cite{Soenjaya}
Let \,$ A \,\subseteq\, X$.\;Then the closure of \,$A$\, is defined as 
\[\overline{A} \,=\, \left\{\, x \,\in\, X \;|\; \,\exists\, \;\{\,x_{\,k}\,\} \,\in\, A \;\;\textit{with}\;  \lim\limits_{k \,\to\, \infty} x_{\,k} \,=\, x \,\right\}.\]
The set \,$ A $\, is said to be closed if $ A \,=\, \overline{A}$. 
\end{definition}

\section{Fixed point theorems in linear $n$-normed space}

\begin{definition}
A sequence \,$\{\,x_{\,k}\,\} \,\subseteq\, X$\, is said to be a b-Cauchy sequence if for every \,$\epsilon \,>\, 0$\, there exists \,$N \,>\, 0$\, such that
\[\text{for every \,$k,\, l \,\geq\, N$},\; \left\|\,x_{\,l} \,-\, x_{\,k},\, b_{\,2},\, \cdots,\, b_{\,n}\,\right\| \,<\, \epsilon.\]
The space \,$X$\, is said to be b-complete if every b-Cauchy sequence is convergent in the semi-normed space \,$\left(\,X,\, \left \|\,\cdot,\, b_{\,2},\, \cdots,\, b_{\,n}\,\right \| \,\right)$.
\end{definition}

\begin{definition}
Let \,$X$\, be a linear \,$n$-normed space and \,$T \,:\, X \,\to\, X$\, be an operator.\,Then the operator \,$T$\, is called b-bounded if there exists some positive constant \,$M$\, such that
\[\left\|\,T\,x,\, b_{\,2},\,\cdots,\,b_{\,n}\,\right\| \,\leq\, M\,\left\|\,x,\, b_{\,2},\, \cdots,\, b_{\,n}\,\right\|\; \;\forall\; x \,\in\, X.\] 
\end{definition}

Let \,$T$\, be a \,$b$-bounded linear operator on \,$X$.\,Then the norm of \,$T$\, is denoted by \,$\|\,T\,\|$\, and is defined as
\begin{align*}
&\|\,T\,\| \,=\, \inf\,\left\{\,M \,:\, \left\|\,T\,x,\, b_{\,2},\,\cdots,\,b_{\,n}\,\right\| \,\leq\, M\,\left\|\,x,\, b_{\,2},\, \cdots,\, b_{\,n}\,\right\|\; \;\forall\; x \,\in\, X\,\right\}.
\end{align*}

\begin{remark}
If \,$T$\, be a \,$b$-bounded on \,$X$, norm of \,$T$\, can be expressed by any one of the following equivalent formula:
\begin{itemize}
\item[(I)]\hspace{.2cm}$\|\,T\,\| \,=\, \sup\,\left\{\,\left\|\,T\,x,\, b_{\,2},\,\cdots,\,b_{\,n}\,\right\| \;:\; \left\|\,x,\, b_{\,2},\, \cdots,\, b_{\,n}\,\right\| \,\leq\, 1\,\right\}$.
\item[(II)]\hspace{.2cm}$\|\,T\,\| \,=\, \sup\,\left\{\,\left\|\,T\,x,\, b_{\,2},\,\cdots,\,b_{\,n}\,\right\| \;:\; \left\|\,x,\, b_{\,2},\, \cdots,\, b_{\,n}\,\right\| \,=\, 1\,\right\}$.
\item[(III)]\hspace{.2cm}$ \|\,T\,\| \,=\, \sup\,\left \{\,\dfrac{\left\|\,T\,x,\, b_{\,2},\,\cdots,\,b_{\,n}\,\right\|}{\left\|\,x,\, b_{\,2},\, \cdots,\, b_{\,n}\,\right\|} \;:\; \left\|\,x,\, b_{\,2},\, \cdots,\, b_{\,n}\,\right\| \,\neq\, 0\,\right \}$. 
\end{itemize}
Also, we have 
\[\left\|\,T\,x,\, b_{\,2},\,\cdots,\,b_{\,n}\,\right\| \,\leq\, \|\,T\,\|\, \left\|\,x,\, b_{\,2},\, \cdots,\, b_{\,n}\,\right\|\, \;\forall\; x \,\in\, X.\]
It is easy to see that \,$X_{b}^{\,\ast}$, collection of all \,$b$-bounded linear operators, forms a Banach space, on \,$X$.
\end{remark}

\begin{definition}
A linear operator \,$T \,:\, X \,\to\, X$\, is said to be continuous at \,$x_{\,0} \,\in\, X$\, if for any open ball \,$B_{\,\{\,e_{\,2},\, \cdots,\, e_{\,n}\,\}}\,(\,T\,(\,x_{\,0}\,),\, \epsilon\,) \,>\, 0$\, there exists a \,$\delta \,>\, 0$\, such that
\[T\,\left(\,B_{\,\{\,e_{\,2},\, \cdots,\, e_{\,n}\,\}}\,(\,x_{\,0},\, \delta\,)\,\right) \,\subset\, B_{\,\{\,e_{\,2},\, \cdots,\, e_{\,n}\,\}}\,(\,T\,(\,x_{\,0}\,),\, \epsilon\,).\]
Equivalently, for given \,$\epsilon \,>\, 0$, there exist some \,$e_{\,2},\, \cdots,\, e_{\,n} \,\in\, X$\, and \,$\delta \,>\, 0$\, such that for \,$x \,\in\, X$
\[\left\|\,x \,-\, x_{\,0},\,e_{\,2},\, \cdots,\, e_{\,n}\,\right\| \,<\, \delta \,\Rightarrow\, \left\|\,T\,(\,x\,) \,-\, T\,(\,x_{\,0}\,),\,e_{\,2},\, \cdots,\, e_{\,n}\,\right\| \,<\, \epsilon.\]
For particular case \,$e_{\,2} \,=\, b_{\,2},\, \cdots,\, e_{\,n} \,=\, b_{\,n}$, it is said to be b-continuous.  
\end{definition}

\begin{definition}
A linear operator \,$T \,:\, X \,\to\, X$\, is said to be b-sequentially continuous at \,$x \,\in\, X$\, if for every sequence \,$\left\{\,x_{\,k}\,\right\}$\, b-converging to \,$x$\, in \,$X$, the sequence \,$\left\{\,T\,(\,x_{\,k}\,)\,\right\}$\, b-converges to \,$T\,(\,x\,)$\, in \,$X$.
\end{definition}

\begin{theorem}
Let \,$T \,:\, X \,\to\, X$\, be a linear operator, where \,$X$\, is a linear n-normed space.\,Then \,$T$\, is b-continuous at \,$x \,\in\, X$\, if and only if \,$T$\, is b-bounded. 
\end{theorem}

\begin{proof}
First we suppose that \,$T$\, is \,$b$-continuous.\,If possible that \,$T$\, is not \,$b$-bounded. Then there exists a sequence \,$\left\{\,x_{\,k}\,\right\}$\, in \,$X$\, such that
\[\left\|\,T\,(\,x_{\,k}\,),\, b_{\,2},\,\cdots,\,b_{\,n}\,\right\| \,>\, k\,\left\|\,x,\, b_{\,2},\, \cdots,\, b_{\,n}\,\right\|\; \;\text{for}\; k \,=\, 1,\,2,\,\cdots.\]
Clearly, \,$x \,\neq\, \theta$, for any \,$k$.\,Let \,$x^{\,\prime}_{\,k} \,=\, \dfrac{x_{\,k}}{\,k\,\left\|\,x_{\,k},\, b_{\,2},\, \cdots,\, b_{\,n}\,\right\|}$.\,Then \,$\left\|\,x^{\,\prime}_{\,k},\, b_{\,2},\, \cdots,\, b_{\,n}\,\right\| \,=\, \dfrac{1}{k} \,\to\, 0$\, as \,$k \,\to\, \infty$, so \,$\left\{\,x^{\,\prime}_{\,k}\,\right\}$\, \,$b$-converges to \,$\theta$.\,Since \,$T$\, is \,$b$-continuous at \,$x \,=\, \theta$, \,$\left\{\,T\,(\,x_{\,k}\,)\,\right\}$\, \,$b$-converges to \,$T\,(\,\theta\,) \,=\, \theta$, i\,.\,e., \,$\left\|\,T\,(\,x^{\,\prime}_{\,k}\,),\, b_{\,2},\,\cdots,\,b_{\,n}\,\right\| \,\to\, 0$\, as \,$k \,\to\, \infty$.\,But on the other hand
\[\left\|\,T\,(\,x^{\,\prime}_{\,k}\,),\, b_{\,2},\,\cdots,\,b_{\,n}\,\right\| \,=\, \dfrac{1}{\,k\,\left\|\,x_{\,k},\, b_{\,2},\, \cdots,\, b_{\,n}\,\right\|}\,\left\|\,T\,(\,x_{\,k}\,),\, b_{\,2},\,\cdots,\,b_{\,n}\,\right\| \,>\, 1,\]for \,$k \,=\, 1,\, 2,\, \cdots$, which is a contradiction.\,Therefore \,$T$\, must be \,$b$-bounded.

Conversely, suppose that \,$T$\, is \,$b$-bounded.\,So, there exists a constant \,$M \,>\, 0$\, such that
\[\left\|\,T\,x,\, b_{\,2},\,\cdots,\,b_{\,n}\,\right\| \,\leq\, M\,\left\|\,x,\, b_{\,2},\, \cdots,\, b_{\,n}\,\right\|\; \;\forall\; x \,\in\, X.\] 
Let \,$\left\{\,x_{\,k}\,\right\}$\, \,$b$-converges to \,$x$, i\,.\,e.,  \,$\left\|\,x_{\,k},\, b_{\,2},\, \cdots,\, b_{\,n}\,\right\| \,\to\, 0$\, as \,$k \,\to\, \infty$.\,Then
\[\left\|\,T\,(\,x_{\,k}\,) \,-\, T\,(\,x\,),\, b_{\,2},\,\cdots,\,b_{\,n}\,\right\| \,=\, \left\|\,T\,\left(\,x_{\,k} \,-\, \,x\,\right),\, b_{\,2},\,\cdots,\,b_{\,n}\,\right\|\]
\[\leq\, M\,\left\|\,x_{\,k} \,-\, x,\, b_{\,2},\, \cdots,\, b_{\,n}\,\right\|\; \;\to\, 0\; \;\text{as \,$k \,\to\, \infty$}.\]
Therefore \,$\left\{\,T\,(\,x_{\,k}\,)\,\right\}$\, \,$b$-converges to \,$T\,(\,x\,)$\, and so \,$T$\, is \,$b$-continuous at \,$x \,\in\, X$.\,Because \,$x \,\in\, X$\, is arbitrary, \,$T$\, is \,$b$-continuous on \,$X$.\,This completes the proof.       
\end{proof}

\begin{theorem}
Let \,$X$\, be a linear n-normed space and \,$T$\, be a linear operator on \,$X$.\;Then \,$T$\, is \,$b$-bounded if and only if \,$T$\, maps bounded sets in \,$X$\, into bounded sets in \,$X$.
\end{theorem}

\begin{proof}
Suppose \,$T$\, is \,$b$-bounded and \,$S$\, is any bounded subset of \,$X$.\;Then there exists \,$M_{\,1} \,>\, 0$\, such that
\[\left\|\,T\,x\,,\, b_{\,2}\,,\, \cdots\,,\, b_{\,n}\,\right\| \,\leq\, M_{\,1}\, \left\|\,x  \,,\, b_{\,2} \,,\, \cdots \,,\, b_{\,n}\,\right\|\; \;\forall\; x \,\in\, X\]and in particular, for all \,$x \,\in\, S$.\;The set \,$S$\, being bounded, for some real number \,$M \,>\, 0$, we have 
\[\left\|\,T\,x\,,\, b_{\,2}\,,\, \cdots\,,\, b_{\,n}\,\right\| \,\leq\, M\; \;\forall\; x \,\in\, S \,\Rightarrow\; \;\text{the set}\; \left\{\,T\,x \,:\, x \,\in\, S \,\right\}\]is bounded in \,$Y$\, and hence \,$T$\, maps bounded sets in \,$X$\, into bounded sets in \,$Y$.

Conversely, for the closed unit ball 
\[B_{\,\{\,e_{\,2} \,,\, \cdots \,,\, e_{\,n}\,\}}\,[\,0 \,,\, 1\,] \,=\, \left\{\,x \,\in\, X \,:\, \left\|\,x \,,\, e_{\,2} \,,\, \cdots \,,\, e_{\,n}\,\right\| \,\leq\, 1\,\right\}\;,\]the set \,$\left\{\,T\,x \,:\, x \,\in\, B_{\,\{\,e_{\,2} \,,\, \cdots \,,\, e_{\,n}\,\}}\,[\,0 \,,\, 1\,]\,\right\}$\, is bounded set in \,$Y$.\,Therefore, there exists \,$K \,>\, 0$\, such that
\[\left\|\,T\,x\,,\, b_{\,2}\,,\, \cdots\,,\, b_{\,n}\,\right\| \,\leq\, K\; \;\;\forall\; x \,\in\, B_{\,\{\,e_{\,2} \,,\, \cdots \,,\, e_{\,n}\,\}}\,[\,0 \,,\, 1\,].\]
If \,$x \,=\, 0$, then \,$T\,x \,=\, 0$\, and the assertion 
\[\left\|\,T\,x\,,\, b_{\,2}\,,\, \cdots\,,\, b_{\,n}\,\right\| \,\leq\, K\, \left\|\,x,\, b_{\,2},\, \cdots,\, b_{\,n}\,\right\|\; \;\text{is obviouly true.\;If}\; x \,\neq\, 0, \,\text{then}\] \,$\dfrac{x}{\|\,x,\, e_{\,2},\, \cdots,\, e_{\,n}\,\|} \,\in\, B_{\{\,e_{\,2},\, \cdots,\, e_{\,n}\,\}}\,[\,0 \,,\, 1\,]$, and for particular \,$e_{\,2} \,=\, b_{\,2},\, \,\cdots,\, e_{\,n} \,=\, b_{\,n}$
\begin{align*}
&\left\|\,T\,\left(\, \dfrac{x}{\|\,x \,,\, b_{\,2} \,,\, \cdots \,,\, b_{\,n}\,\|}\,\right)\,\right\| \,\leq\, K\\
&\Rightarrow\, \left\|\,T\,x \,,\, b_{\,2}\,,\, \cdots\,,\, b_{\,n}\,\right\| \,\leq\, K\, \left\|\,x  \,,\, b_{\,2} \,,\, \cdots \,,\, b_{\,n}\,\right\|\; \;\forall\; x \,\in\, X.
\end{align*}
Hence, \,$T$\, is a \,$b$-bounded linear operator.   
\end{proof} 

\begin{definition}
Let \,$X$\, be a linear \,$n$-normed space and \,$T \,:\, X \,\to\, X$\, be an operator.\,The operator \,$T$\, is called b-contraction operator if
\begin{equation}\label{eq1}
\left\|\,T\,x \,-\, T\,y,\, b_{\,2},\,\cdots,\,b_{\,n}\,\right\| \,\leq\, \alpha\,\left\|\,x \,-\, y,\, b_{\,2},\, \cdots,\, b_{\,n}\,\right\|\; \;\forall\; x,\,y \,\in\, X,
\end{equation}
where \,$0 \,<\, \alpha \,<\, 1$.
\end{definition}

\begin{theorem}
Any b-contraction operator is b-continuous.
\end{theorem}

\begin{proof}
Let \,$T \,:\, X \,\to\, X$\, be a \,$b$-contraction operator.\,Choose \,$\epsilon \,>\, 0$\, arbitrary and take \,$0 \,<\, \delta \,<\, \dfrac{\epsilon}{\alpha}$.\,Then for \,$x,\, y \,\in\, X$, we have
\begin{align*}
&\left\|\,T\,x \,-\, T\,y,\, b_{\,2},\,\cdots,\,b_{\,n}\,\right\| \,\leq\, \alpha\,\left\|\,x \,-\, y,\, b_{\,2},\, \cdots,\, b_{\,n}\,\right\|\\
& \,<\, \alpha\,\delta\;\; \;\text{if}\;\;\left\|\,x \,-\, y,\, b_{\,2},\, \cdots,\, b_{\,n}\,\right\| \,<\, \delta.
\end{align*} 
So, \,$\left\|\,T\,x \,-\, T\,y,\, b_{\,2},\,\cdots,\,b_{\,n}\,\right\| \,<\, \epsilon$\, whenever \,$\left\|\,x \,-\, y,\, b_{\,2},\, \cdots,\, b_{\,n}\,\right\| \,<\, \delta,\; x,\,y \,\in\, X$.\,Therefore, \,$T$\, is \,$b$-continuous. 
\end{proof}

\begin{theorem}\label{th1}
Let \,$X$\, be a \,$b$-complete linear \,$n$-normed space and \,$T \,:\, X \,\to\, X$\, be an \,$b$-contraction mapping.\,Then there exists a unique fixed point \,$x_{\,0}$\, of the operator \,$T$\, in \,$X$\, provided the set \,$\left\{\,x_{\,0},\, b_{\,2},\, \cdots,\, b_{\,n}\,\right\}$\, is linearly independent.
\end{theorem}

\begin{proof}
Let \,$x \,\in\, X$\, be an arbitrary element.\;Starting from \,$x$, we form the iterations
\[x_{\,1} \,=\, T\,x,\, x_{\,2} \,=\, T\,x_{\,1},\, x_{\,3} \,=\, T\,x_{\,2},\, \cdots,\, ,\, x_{\,k} \,=\, T\,x_{\,k \,-\, 1},\, \cdots. \]
We verify that \,$\left\{\,x_{\,k}\,\right\}$\, is a \,$b$-Cauchy sequence.\,We have
\begin{align*}
\left\|\,x_{\,1} \,-\, x_{\,2},\, b_{\,2},\, \cdots,\, b_{\,n}\,\right\| &\,=\, \left\|\,T\,x \,-\, T\,x_{\,1},\, b_{\,2},\,\cdots,\,b_{\,n}\,\right\|\\
& \,\leq\, \alpha\,\left\|\,x \,-\, x_{\,1},\, b_{\,2},\, \cdots,\, b_{\,n}\,\right\|\\
& =\, \alpha\,\left\|\,x \,-\, T\,x,\, b_{\,2},\, \cdots,\, b_{\,n}\,\right\|.
\end{align*}
\begin{align*}
\left\|\,x_{\,2} \,-\, x_{\,3},\, b_{\,2},\, \cdots,\, b_{\,n}\,\right\|& \,=\, \left\|\,T\,x_{\,1} \,-\, T\,x_{\,2},\, b_{\,2},\,\cdots,\,b_{\,n}\,\right\|\\
& \,\leq\, \alpha\,\left\|\,x_{\,1} \,-\, x_{\,2},\, b_{\,2},\, \cdots,\, b_{\,n}\,\right\|\\
&=\, \alpha^{\,2}\,\left\|\,x \,-\, T\,x,\, b_{\,2},\, \cdots,\, b_{\,n}\,\right\|.
\end{align*}
\begin{align*}
\left\|\,x_{\,3} \,-\, x_{\,4},\, b_{\,2},\, \cdots,\, b_{\,n}\,\right\| &\,=\, \left\|\,T\,x_{\,2} \,-\, T\,x_{\,3},\, b_{\,2},\,\cdots,\,b_{\,n}\,\right\|\\
& \,\leq\, \alpha\,\left\|\,x_{\,2} \,-\, x_{\,3},\, b_{\,2},\, \cdots,\, b_{\,n}\,\right\|\\
&=\, \alpha^{\,3}\,\left\|\,x \,-\, T\,x,\, b_{\,2},\, \cdots,\, b_{\,n}\,\right\|.
\end{align*}
In general, for any positive integer \,$k$,
\[\left\|\,x_{\,k} \,-\, x_{\,k \,+\, 1},\, b_{\,2},\, \cdots,\, b_{\,n}\,\right\| \,\leq\, \alpha^{\,k}\,\left\|\,x \,-\, T\,x,\, b_{\,2},\, \cdots,\, b_{\,n}\,\right\|.\] 
Also, for any positive integer \,$p$, we have
\begin{align}
&\left\|\,x_{\,k} \,-\, x_{\,k \,+\, p},\, b_{\,2},\, \cdots,\, b_{\,n}\,\right\|\nonumber\\
&\leq\,\left\|\,x_{\,k} - x_{\,k \,+\, 1},\, b_{\,2},\, \cdots,\, b_{\,n}\,\right\| \,+\, \left\|\,x_{\,k \,+\, 1} - x_{\,k \,+\, 2},\, b_{\,2},\, \cdots,\, b_{\,n}\,\right\|\,+\,\cdots\nonumber\\
&\hspace{2cm}\cdots\,+\,\left\|\,x_{\,k\,+\,p\,-1} - x_{\,k \,+\, p},\, b_{\,2},\, \cdots,\, b_{\,n}\,\right\|\nonumber\\
&\leq\, \alpha^{\,k}\,\left\|\,x \,-\, T\,x,\, b_{\,2},\, \cdots,\, b_{\,n}\,\right\| \,+\, \alpha^{\,k\,+\,1}\,\left\|\,x \,-\, T\,x,\, b_{\,2},\, \cdots,\, b_{\,n}\,\right\|\,+\,\cdots\nonumber\\
&\hspace{2cm}\cdots\,+\,\alpha^{\,k\,+\,p\,-\,1}\,\left\|\,x \,-\, T\,x,\, b_{\,2},\, \cdots,\, b_{\,n}\,\right\|\nonumber\\
&=\, \left(\,\alpha^{\,k} \,+\, \alpha^{\,k\,+\,1} \,+\, \cdots\,+\, \alpha^{\,k\,+\,p\,-\,1}\,\right)\,\left\|\,x \,-\, T\,x,\, b_{\,2},\, \cdots,\, b_{\,n}\,\right\|\nonumber\\
&=\,\dfrac{\alpha^{\,k} \,-\, \alpha^{\,k\,+\,p}}{1 \,-\, \alpha}\,\left\|\,x \,-\, T\,x,\, b_{\,2},\, \cdots,\, b_{\,n}\,\right\|\label{eq1.1}\\
&<\,\dfrac{\alpha^{\,k}}{1 \,-\, \alpha}\,\left\|\,x \,-\, T\,x,\, b_{\,2},\, \cdots,\, b_{\,n}\,\right\|,\;[\;\text{since}\; \,0 \,<\, \alpha \,<\, 1\;].\label{eq1.2}
\end{align}
Since \,$\alpha \,<\, 1$, the relation (\ref{eq1.2}) shows that \,$\left\|\,x_{\,k} \,-\, x_{\,k \,+\, p},\, b_{\,2},\, \cdots,\, b_{\,n}\,\right\| \,\to\, 0$\, as \,$k \,\to\, \infty$.\;Therefore the sequence \,$\left\{\,x_{\,k}\,\right\}$\, is a \,$b$-Cauchy sequence.\,Since \,$X$\, is \,$b$-complete linear \,$n$-normed space, the sequence \,$\left\{\,x_{\,k}\,\right\}$\, is convergent in the semi-normed space \,$\left(\,X,\, \left\|\,\cdot,\,b_{\,2},\, \cdots,\, b_{\,n}\,\right\|\,\right)$.\;So, let \,$\lim\limits_{k \,\to\, \infty}\,x_{\,k} \,=\, x_{\,0}$\, with the property that the set \,$\left\{\,x_{\,0},\, b_{\,2},\, \cdots,\, b_{\,n}\,\right\}$\, is linearly independent.\,We now show that \,$x_{\,0}$\, is a fixed point of the operator \,$T$, i\,.\,e., \,$T\,x_{\,0} \,=\, x_{\,0}$.\;We have
\begin{align*}
&\left\|\,x_{\,0} \,-\, T\,x_{\,0},\, b_{\,2},\, \cdots,\, b_{\,n}\,\right\| \,\leq\, \left\|\,x_{\,0} \,-\, x_{\,k},\, b_{\,2},\, \cdots,\, b_{\,n}\,\right\| \,+\, \left\|\,x_{\,k} \,-\, T\,x_{\,0},\, b_{\,2},\, \cdots,\, b_{\,n}\,\right\|\\
&=\, \left\|\,x_{\,0} \,-\, x_{\,k},\, b_{\,2},\, \cdots,\, b_{\,n}\,\right\| \,+\, \left\|\,T\,x_{\,k\,-\,1} \,-\, T\,x_{\,0},\, b_{\,2},\, \cdots,\, b_{\,n}\,\right\|\\
&\leq\, \left\|\,x_{\,0} \,-\, x_{\,k},\, b_{\,2},\, \cdots,\, b_{\,n}\,\right\| \,+\, \alpha\,\left\|\,x_{\,k\,-\,1} \,-\, x_{\,0},\, b_{\,2},\, \cdots,\, b_{\,n}\,\right\|\; [\;\text{by (\ref{eq1})}\;]\\
&\to\, 0\; \;\text{as}\; \,k \,\to\, \infty\; \;\left[\;\text{since}\;\lim\limits_{k \,\to\, \infty}\,x_{\,k} \,=\, x_{\,0}\,\right].
\end{align*}
So, \,$T\,x_{\,0} \,=\, x_{\,0}$.\;Therefore \,$x_{\,0}$\, is a fixed point of \,$T$.\;We now verify that there exists only one fixed point of \,$T$.\,Let \,$y_{\,0} \,\in\, X$\, with \,$\left\{\,x_{\,0},\, b_{\,2},\, \cdots,\, b_{\,n}\,\right\}$\, is linearly independent such that \,$T\,y_{\,0} \,=\, y_{\,0}$.\,Then using (\ref{eq1}), we have
\begin{align*}
\left\|\,x_{\,0} \,-\, y_{\,0},\, b_{\,2},\, \cdots,\, b_{\,n}\,\right\| &\,=\, \left\|\,T\,x_{\,0} \,-\, T\,y_{\,0},\, b_{\,2},\, \cdots,\, b_{\,n}\,\right\|\\
& \,\leq\, \alpha\,\left\|\,x_{\,0} \,-\, y_{\,0},\, b_{\,2},\, \cdots,\, b_{\,n}\,\right\|.
\end{align*}
If \,$\left\|\,x_{\,0} \,-\, y_{\,0},\, b_{\,2},\, \cdots,\, b_{\,n}\,\right\| \,>\, 0$\, then from above inequality, we obtain \,$\alpha \,\geq\, 1$, which is a contradiction.\,Hence, \,$\left\|\,x_{\,0} \,-\, y_{\,0},\, b_{\,2},\, \cdots,\, b_{\,n}\,\right\| \,=\, 0$, i\,.\,e., \,$x_{\,0} \,=\, y_{\,0}$\, and so \,$T$\, has a unique fixed point in \,$X$.\,This proves the theorem.     
\end{proof}

\begin{note}
We consider the inequality (\ref{eq1.1})
\[\left\|\,x_{\,k} \,-\, x_{\,k \,+\, p},\, b_{\,2},\, \cdots,\, b_{\,n}\,\right\| \,\leq\, \dfrac{\alpha^{\,k} \,-\, \alpha^{\,k\,+\,p}}{1 \,-\, \alpha}\,\left\|\,x \,-\, T\,x,\, b_{\,2},\, \cdots,\, b_{\,n}\,\right\|.\]
As \,$p \,\to\, \infty$, since \,$\alpha \,<\, 1$, the right hand side tends to \,$\dfrac{\alpha^{\,k}}{1 \,-\, \alpha}\,\left\|\,x \,-\, T\,x,\, b_{\,2},\, \cdots,\, b_{\,n}\,\right\|$\, and the left hand side tends to \,$\left\|\,x_{\,k} \,-\, x_{\,0},\, b_{\,2},\, \cdots,\, b_{\,n}\,\right\|$\, because \,$x_{\,k \,+\, p} \,\to\, x_{\,0}$.\,So,
\begin{equation}\label{eqn1.3}
\left\|\,x_{\,k} \,-\, x_{\,0},\, b_{\,2},\, \cdots,\, b_{\,n}\,\right\| \,\leq\, \dfrac{\alpha^{\,k}}{1 \,-\, \alpha}\,\left\|\,x \,-\, T\,x,\, b_{\,2},\, \cdots,\, b_{\,n}\,\right\|.
\end{equation}
The relation (\ref{eqn1.3}) gives an estimation for the error of the \,$k\text{th}$\, approximation.
\end{note}

Sometimes it may happen that although the condition (\ref{eq1}) is not satisfied on the entire space \,$X$, it is satisfied on a certain subset of \,$X$.\,Under such situation, we prove the following theorem. 

\begin{theorem}
Let \,$X$\, be a \,$b$-complete linear \,$n$-normed space and \,$T \,:\, X \,\to\, X$\, be an operator such that
\[\left\|\,T\,x \,-\, T\,y,\, b_{\,2},\,\cdots,\,b_{\,n}\,\right\| \,\leq\, \alpha\,\left\|\,x \,-\, y,\, b_{\,2},\, \cdots,\, b_{\,n}\,\right\|,\]
for \,$x,\, y \,\in\, \overline{\,B} \,=\, B_{\,\{\,b_{\,2},\, \cdots,\, b_{\,n}\,\}}\,[\,x_{\,0},\,r\,]$.\,Moreover, assume that \[\left\|\,x_{\,0} \,-\, T\,x_{\,0},\, b_{\,2},\, \cdots,\, b_{\,n}\,\right\| \,<\, (\,1 \,-\, \alpha\,)\,r.\] Then the iterative sequence starting from \,$x_{\,0}$\, converges to an \,$x \,\in\, \overline{\,B}$\, which is a unique fixed point of \,$T$\, in \,$\overline{\,B}$  
\end{theorem}

\begin{proof}
We verify by induction that \,$x_{\,1} \,=\, T\,x_{\,0},\, x_{\,2} \,=\, T\,x_{\,1},\, x_{\,3} \,=\, T\,x_{\,2},\, \cdots,\, x_{\,k} \,=\, T\,x_{\,k \,-\, 1},\, \cdots$\, are in \,$\overline{\,B}$.\,Clearly, \,$x_{\,0} \,\in\, \overline{\,B}$.\,Assume that \,$x_{\,1},\, x_{\,2},\, \cdots,\, x_{\,k \,-\, 1}$\, are in \,$\overline{\,B}$.\,We show that \,$x_{\,k}$\, also lies in \,$\overline{\,B}$.\,We have
\begin{align*}
&\left\|\,x_{\,2} \,-\, x_{\,1},\, b_{\,2},\, \cdots,\, b_{\,n}\,\right\| \,\leq\, \alpha\,\left\|\,x_{\,1} \,-\, x_{\,0},\, b_{\,2},\, \cdots,\, b_{\,n}\,\right\| \,\leq\, \alpha\,(\,1 \,-\, \alpha\,)\,r\\
&\left\|\,x_{\,3} \,-\, x_{\,2},\, b_{\,2},\, \cdots,\, b_{\,n}\,\right\| \,\leq\, \alpha\,\left\|\,x_{\,2} \,-\, x_{\,1},\, b_{\,2},\, \cdots,\, b_{\,n}\,\right\| \,\leq\, \alpha^{\,2}\,(\,1 \,-\, \alpha\,)\,r\\
&\hspace{.3cm}\cdots\hspace{6.5cm}\cdots\hspace{4cm}\cdots\\
&\left\|\,x_{\,k\,-\,1} \,-\, x_{\,k\,-\,2},\, b_{\,2},\, \cdots,\, b_{\,n}\,\right\| \,\leq\, \alpha\,\left\|\,x_{\,k\,-\,2} \,-\, x_{\,k\,-\,3},\, b_{\,2},\, \cdots,\, b_{\,n}\,\right\| \,\leq\, \alpha^{\,k\,-\,2}\,(\,1 \,-\, \alpha\,)\,r\\
&\left\|\,x_{\,k} \,-\, x_{\,k\,-\,1},\, b_{\,2},\, \cdots,\, b_{\,n}\,\right\| \,\leq\, \alpha\,\left\|\,x_{\,k\,-\,1} \,-\, x_{\,k\,-\,2},\, b_{\,2},\, \cdots,\, b_{\,n}\,\right\| \,\leq\, \alpha^{\,k\,-\,1}\,(\,1 \,-\, \alpha\,)\,r.
\end{align*}
Therefore, 
\begin{align*}
&\left\|\,x_{\,0} \,-\, x_{\,k},\, b_{\,2},\, \cdots,\, b_{\,n}\,\right\|\\
&\leq\, \left\|\,x_{\,0} - x_{\,1},\, b_{\,2},\, \cdots,\, b_{\,n}\,\right\| \,+\, \left\|\,x_{\,1} - x_{\,2},\, b_{\,2},\, \cdots,\, b_{\,n}\,\right\|\, + \cdots + \left\|\,x_{\,k\,-\,1} - x_{\,k},\, b_{\,2},\, \cdots,\, b_{\,n}\,\right\|\\
&\leq\,(\,1 \,-\, \alpha\,)\,r \,+\, \alpha\,(\,1 \,-\, \alpha\,)\,r \,+\, \alpha^{\,2}\,(\,1 \,-\, \alpha\,)\,r \,+\, \cdots \,+\,\alpha^{\,k\,-\,1}\,(\,1 \,-\, \alpha\,)\,r\\
&=\,(\,1 \,-\, \alpha\,)\,r\left(\,1 \,+\, \alpha \,+\, \alpha^{\,2} \,+\, \cdots \,+\, \alpha^{\,k\,-\,1}\,\right) \,=\, \left(\,1 \,-\, \alpha^{\,k}\,\right)\,r \,<\, r.
\end{align*}
So, \,$x_{\,k}$\, lies in \,$\overline{\,B}$.\,Therefore, every member of the sequence \,$\left\{\,x_{\,k}\,\right\}$\, is contained in \,$\overline{\,B}$.\,According to the Theorem \ref{th1}, we now see that the sequence \,$\left\{\,x_{\,k}\,\right\}$\, converges to an element \,$x^{\,\prime}_{\,0}$, say.\,The element \,$x^{\,\prime}_{\,0}$\, belong to \,$\overline{\,B}$\, because \,$\overline{\,B}$\, is closed.\,The method of proof of Theorem \ref{th1}, shows that \,$x^{\,\prime}_{\,0}$\, is the unique fixed point of \,$T$\, in \,$\overline{\,B}$.\,This completes the proof of the theorem. 
\end{proof}

The Cartesian product \,$X \,\times\, X$\, of linear\;$n$-normed space \,$X$\, is a linear\;$n$-normed space with respect to the \,$n$-norm given by
\[ \left\|\,(\,x_{\,1},\, y_{\,1}\,),\, (\,x_{\,2},\, y_{\,2}\,),\, \cdots,\, (\,x_{\,n},\, y_{\,n}\,)\,\right\|_{1} \,=\, \,\|\,x_{\,1},\, x_{\,2},\, \cdots,\, x_{\,n}\,\| \,+\, \|\,y_{\,1},\, y_{\,2},\, \cdots,\, y_{\,n}\,\|,\]
for all \,$(\,x_{\,1},\, y_{\,1}\,),\, (\,x_{\,2},\, y_{\,2}\,),\, \cdots,\, (\,x_{\,n},\, y_{\,n}\,) \,\in\, X \,\times\, X$.

\begin{lemma}\label{lem1}
Let \,$X$\, be a linear \,$n$-normed space and \,$X \,\times\, X$\, be the product linear \,$n$-normed space with respect to the n-norm \,$\left \|\,\cdot,\, \cdots,\, \cdot \,\right \|_{1}$.\,Consider the open ball \,$B_{\,\{\,a_{\,2},\, \cdots,\, a_{\,n}\,\}}\,\left(\,(\,x_{\,0},\, y_{\,0}\,),\, r_{\,1}\,\right)$\, in \,$X \,\times\, X$, where \,$a_{\,2} \,=\, \left(\,e_{\,2},\,e^{\,\prime}_{\,2}\,\right),\, \cdots,\, a_{\,n} \,=\, \left(\,e_{\,n},\,e^{\,\prime}_{\,n}\,\right)$.\,Then there exist two open balls \,$B_{\,\{\,e_{\,2},\, \cdots,\, e_{\,n}\,\}}\,(\,x_{\,0},\, r\,)$\, and \,$B_{\,\{\,e^{\,\prime}_{\,2},\, \cdots,\, e^{\,\prime}_{\,n}\,\}}\,(\,y_{\,0},\, r^{\,\prime}\,)$\, in \,$X$\, such that
\[B_{\,\{\,e_{\,2},\, \cdots,\, e_{\,n}\,\}}\,(\,x_{\,0},\, r\,) \,\times\, B_{\,\{\,e^{\,\prime}_{\,2},\, \cdots,\, e^{\,\prime}_{\,n}\,\}}\,(\,y_{\,0},\, r^{\,\prime}\,) \,\subseteq\, B_{\,\{\,a_{\,2},\, \cdots,\, a_{\,n}\,\}}\,\left(\,(\,x_{\,0},\, y_{\,0}\,),\, r_{\,1}\,\right)\]     
\end{lemma}

\begin{proof}
Select \,$r$\, and \,$r^{\,\prime}$\, positive numbers both less than  \,$r_{\,1}$.\,Let 
\[x \,\in\, B_{\,\{\,e_{\,2},\, \cdots,\, e_{\,n}\,\}}\,(\,x_{\,0},\, r\,)\; \;\text{and}\; \,y \,\in\, B_{\,\{\,e^{\,\prime}_{\,2},\, \cdots,\, e^{\,\prime}_{\,n}\,\}}\,(\,y_{\,0},\, r^{\,\prime}\,).\]
Then
\[\left\|\,x \,-\, x_{\,0},\, e_{\,2},\, \cdots,\, e_{\,n}\,\right\| \,<\, r\; \;\text{and}\; \;\left\|\,y \,-\, y_{\,0},\, e^{\,\prime}_{\,2},\, \cdots,\, e^{\,\prime}_{\,n}\,\right\| \,<\, r^{\,\prime}.\]Now
\begin{align*}
&\left\|\,(\,x,\,y\,) \,-\, (\,x_{\,0},\,y_{\,0}\,),\, (\,e_{\,2},\,e^{\,\prime}_{\,2}\,),\, \cdots,\, (\,e_{\,n},\,e^{\,\prime}_{\,n}\,)\,\right\|\\
& \,=\, \left\|\,(\,x \,-\, x_{\,0}\,),\, (\,y \,-\, y_{\,0}\,),\, (\,e_{\,2},\,e^{\,\prime}_{\,2}\,),\, \cdots,\, (\,e_{\,n},\,e^{\,\prime}_{\,n}\,)\,\right\|\\
&=\, \left\|\,x \,-\, x_{\,0},\, e_{\,2},\, \cdots,\, e_{\,n}\,\right\| \,+\, \left\|\,y \,-\, y_{\,0},\, e^{\,\prime}_{\,2},\, \cdots,\, e^{\,\prime}_{\,n}\,\right\| \,<\, r \,+\, r^{\,\prime} \,<\, r_{\,1}.
\end{align*}
So, \,$(\,x,\,y\,) \,\in\, B_{\,\{\,a_{\,2},\, \cdots,\, a_{\,n}\,\}}\,\left(\,(\,x_{\,0},\, y_{\,0}\,),\, r_{\,1}\,\right)$.\,This proves the lemma. 
\end{proof}

\begin{theorem}
Let \,$T$\, be a linear operator on a \,$X$\, such that
\begin{equation}\label{eq1.4}
\left\|\,T\,x \,-\, T\,y,\, b_{\,2},\,\cdots,\,b_{\,n}\,\right\| \,\leq\, \left\|\,x \,-\, y,\, b_{\,2},\, \cdots,\, b_{\,n}\,\right\|\; \;\forall\; x,\,y \,\in\, X, \,x \,\neq\, y.
\end{equation}
Suppose that there exists a point \,$x \,\in\, X$\, such that the sequence of iterates \,$\left\{\,T^{\,k}\,(\,x\,)\,\right\}$\, has a subsequence \,$b$-converging to \,$\xi \,\in\, X$.\,Then \,$\xi$\, is the unique fixed point of \,$T$. 
\end{theorem}

\begin{proof}
Let \,$\left\{\,T^{\,k_{\,i}}\,(\,x\,)\,\right\}$\, be that subsequence of \,$\left\{\,T^{\,k}\,(\,x\,)\,\right\}$\, which is \,$b$-convergent.\,So, take 
\begin{equation}\label{eq1.5}
\lim\limits_{i \,\to\, \infty}\,T^{\,k_{\,i}}\,(\,x\,) \,=\, \xi\; \;\text{(\,say\,)}.
\end{equation}
Suppose, if possible that \,$\xi$\, is not a fixed point of \,$T$, i\,.\,e., \,$T\,(\,\xi\,) \,\neq\, \xi$.\,Let \,$Y$\, be the subset of \,$X \,\times\, X$\, defined by 
\[Y \,=\, X \,\times\, X \,-\, \Delta,\; \;\text{where}\; \;\Delta \,=\, \left\{\,(\,x,\,y\,) \,\in\, X \,\times\, X \,:\, x \,=\, y\,\right\}.\]
We now define a real-valued function of two variables on \,$Y$\, by
\[f\,(\,p,\,q\,) \,=\, \dfrac{\left\|\,T\,(\,p\,) \,-\, T\,(\,q\,),\, b_{\,2},\,\cdots,\,b_{\,n}\,\right\|}{\left\|\,p \,-\, q,\, b_{\,2},\,\cdots,\,b_{\,n}\,\right\|},\; \;(\,p,\,q\,) \,\in\, Y.\]
If \,$\left\{\,p_{\,k}\,\right\}$\, \,$b$-converges to \,$p$\, and \,$\left\{\,q_{\,k}\,\right\}$\, \,$b$-converges to \,$q$, then because \,$T$\, is \,$b$-continuous, we get that \,$\left\{\,T\,(\,p_{\,k}\,)\,\right\}$\, \,$b$-converges to \,$T\,(\,p\,)$\, and \,$\left\{\,T\,(\,q_{\,k}\,)\,\right\}$\, \,$b$-converges to \,$T\,(\,q\,)$\, and so we have
\begin{align*}
f\,(\,p_{\,k},\,q_{\,k}\,) &\,=\, \dfrac{\left\|\,T\,(\,p_{\,k}\,) \,-\, T\,(\,q_{\,k}\,),\, b_{\,2},\,\cdots,\,b_{\,n}\,\right\|}{\left\|\,p_{\,k} \,-\, q_{\,k},\, b_{\,2},\,\cdots,\,b_{\,n}\,\right\|}\\
& \,\to\, \dfrac{\left\|\,T\,(\,p\,) \,-\, T\,(\,q\,),\, b_{\,2},\,\cdots,\,b_{\,n}\,\right\|}{\left\|\,p \,-\, q,\, b_{\,2},\,\cdots,\,b_{\,n}\,\right\|}
\end{align*} 
for \,$(\,p_{\,k},\,q_{\,k}\,),\, (\,p,\,q\,) \,\in\, Y$.\,This implies that \,$f\,(\,p,\,q\,)$\, is \,$b$-continuous on \,$Y$.\,By (\ref{eq1.4}), \,$f\,(\,p,\,q\,) \,<\, 1$\, for \,$(\,p,\,q\,) \,\in\, Y$\, and so \,$f\left(\,\xi,\,T\,(\,\xi\,)\,\right) \,<\, 1$.\,Since \,$f\,(\,p,\,q\,)$\, is \,$b$-continuous at \,$(\,\xi,\,T\,(\,\xi\,)\,)$, if \,$f\left(\,\xi,\,T\,(\,\xi\,)\,\right) \,<\, R \,<\, 1$, there exists a open ball \,$U$, say arround \,$(\,\xi,\,T\,(\,\xi\,)\,)$\, such that
\begin{equation}\label{eq1.6}
\text{for}\; \;(\,p,\,q\,) \,\in\, U \,\cap\, Y, \;0 \,\leq\, f\,(\,p,\,q\,) \,\leq\, R \,<\, 1.
\end{equation}
By lemma \ref{lem1}, there exists two open balls \,$B_{\,1} \,=\, B_{\,\{\,b_{\,2},\, \cdots,\, b_{\,n}\,\}}\,(\,\xi,\, \alpha\,)$\, and \,$B_{\,2} \,=\, B_{\,\{\,b_{\,2},\, \cdots,\, b_{\,n}\,\}}\,(\,T\,\xi,\, \alpha\,)$\, such that \,$B_{\,1} \,\times\, B_{\,2} \,\subset\, U$.\,The positive number \,$\alpha$\, may be choosen small enough to ensure also \,$\alpha \,<\, \dfrac{1}{3}\,\left\|\,\xi \,-\, T\,\xi,\, b_{\,2},\, \cdots,\, b_{\,n}\,\right\|$\, and this only shows that the open balls \,$B_{1}$\, and \,$B_{2}$\, are disjoint.\,By (\ref{eq1.5}), there exists a positive integer \,$N$\, such that \,$T^{\,k_{\,i}}\,(\,x\,) \,\in\, B_{1}$\, for \,$i \,>\, N$.\,By (\ref{eq1.4})
\[\left\|\,T^{\,k_{\,i} \,+\, 1}\,(\,x\,) \,-\, T\,\xi,\, b_{\,2},\,\cdots,\,b_{\,n}\,\right\| \,<\, \left\|\,T^{\,k_{\,i}}\,(\,x\,) \,-\, \xi,\, b_{\,2},\,\cdots,\,b_{\,n}\,\right\| \,<\, \alpha\]
so that \,$T^{\,k_{\,i} \,+\, 1}\,(\,x\,) \,\in\, B_{2}$\, for \,$i \,>\, N$.\,Since \,$B_{1} \,\cap\, B_{2} \,=\, \phi$, it follows that 
\begin{equation}\label{eq1.7}
\left\|\,T^{\,k_{\,i}}\,(\,x\,) \,-\, T^{\,k_{\,i} \,+\, 1}\,(\,x\,),\, b_{\,2},\,\cdots,\,b_{\,n}\,\right\| \,>\, \alpha\; \;\text{for}\; \;i \,>\, N.
\end{equation}
Now, for \,$i \,>\, N$, \,$T^{\,k_{\,i}}\,(\,x\,) \,\in\, B_{1}$\, and \,$T^{\,k_{\,i} \,+\, 1}\,(\,x\,) \,\in\, B_{2}$\, and so from (\ref{eq1.6}) \\$f\,\left(\,T^{\,k_{\,i}}\,(\,x\,),\,T^{\,k_{\,i} \,+\, 1}\,(\,x\,)\,\right) \,<\, R$, i\,.\,e.,
\begin{align}
&\left\|\,T^{\,k_{\,i} \,+\, 1}\,(\,x\,) \,-\, T^{\,k_{\,i} \,+\, 2}\,(\,x\,),\, b_{\,2},\,\cdots,\,b_{\,n}\,\right\|\nonumber\\
& \,<\, R\,\left\|\,T^{\,k_{\,i}}\,(\,x\,) \,-\, T^{\,k_{\,i} \,+\, 1}\,(\,x\,),\, b_{\,2},\,\cdots,\,b_{\,n}\,\right\|.\label{eq1.81}
\end{align}
Let \,$l \,>\, j \,>\, N$.\,A repeated use of (\ref{eq1.4}) and (\ref{eq1.81}) gives
\begin{align*}
&\left\|\,T^{\,k_{\,l}}\,(\,x\,) \,-\, T^{\,k_{\,l} \,+\, 1}\,(\,x\,),\, b_{\,2},\,\cdots,\,b_{\,n}\,\right\|\\
& \,\leq\, \left\|\,T^{\,k_{\,l \,-\,1} \,+\, 1}\,(\,x\,) \,-\, T^{\,k_{\,l \,-\,l} \,+\, 2}\,(\,x\,),\, b_{\,2},\,\cdots,\,b_{\,n}\,\right\|\\
&<\, R\,\left\|\,T^{\,k_{\,l \,-\,1}}\,(\,x\,) \,-\, T^{\,k_{\,l \,-\,l} \,+\, 1}\,(\,x\,),\, b_{\,2},\,\cdots,\,b_{\,n}\,\right\|\\
&\leq\, \cdots\,\hspace{1.5cm}\,\cdots\,\hspace{1cm}\,\cdots\hspace{1.5cm},\,\cdots\\
&\leq\, R^{\,l \,-\, j}\,\left\|\,T^{\,k_{\,j}}\,(\,x\,) \,-\, T^{\,k_{\,j} \,+\, 1}\,(\,x\,),\, b_{\,2},\,\cdots,\,b_{\,n}\,\right\| \,\to\, 0\; \;\text{as}\; \;k \,\to\, \infty\; \;[\;\text{since}\;R \,<\, 1\;].
\end{align*} 
But this last relation contradicts (\ref{eq1.7}).\,Therefore we arrive at a contradiction.\,Hence, \,$T\,(\,\xi\,) \,=\, \xi$\, and \,$\xi$\, is a fixed point of \,$T$.\,If \,$\eta \,\neq\, \xi$\, is another fixed point i\,.\,e., \,$T\,(\,\eta\,) \,=\, \eta$, then
\[\left\|\,T\,\xi \,-\, T\,\eta,\, b_{\,2},\,\cdots,\,b_{\,n}\,\right\| \,\leq\, \left\|\,\xi \,-\, \eta,\, b_{\,2},\, \cdots,\, b_{\,n}\,\right\|\]
against (\ref{eq1.4}).\,Hence, \,$\xi$\, is the unique fixed point \,$T$.\,This proves the theorem.   
\end{proof}

Next, we prove some further theorems on fixed points of operators in linear \,$n$-normed space.

\begin{theorem}\label{thmm1}
Let \,$X$\, be a \,$b$-complete linear \,$n$-normed space and \,$T$\, be a mapping of \,$X$\, into itself.\,Suppose that for each positive integer \,$k$,
 \begin{equation}\label{eq1.8}
\left\|\,T^{\,k}\,x \,-\, T^{\,k}\,y,\, b_{\,2},\,\cdots,\,b_{\,n}\,\right\| \,\leq\, a_{\,k}\,\left\|\,x \,-\, y,\, b_{\,2},\, \cdots,\, b_{\,n}\,\right\|\; \;\forall\; x,\,y \,\in\, X,
\end{equation}
where \,$a_{\,k} \,>\, 0$\, is independent of \,$x,\,y$.\,If the series \,$\sum\limits_{k \,=\, 1}^{\,\infty}\,a_{\,k}$\, is convergent, then \,$T$\, has a unique fixed point in \,$X$.
\end{theorem}

\begin{proof}
Let \,$x_{\,0}$\, be an arbitrary element in \,$X$\, and consider the sequence \,$\left\{\,x_{\,k}\,\right\}$\, of iterates, \,$x_{\,k} \,=\, T^{\,k}\,x_{\,0},\; \;k \,=\, 1,\,2,\,3,\cdots$.\,We note that \,$x_{k \,+\, 1} \,=\, T^{\,k \,+\, 1}\,x_{\,0} \,=\, T^{\,k}\,\left(\,T\,x_{\,0}\,\right) \,=\, T^{\,k}\,x_{\,1}$\, and also, \,$x_{k \,+\, 1} \,=\, T\,\left(\,T^{\,k}\,x_{\,0}\,\right) \,=\, T\,x_{\,k}$.\,If \,$p$\, and \,$q\;(\,p \,>\, q\,)$\, be positive integers, then we obtain
\begin{align*}
\left\|\,x_{\,p} \,-\, x_{\,q},\, b_{\,2},\, \cdots,\, b_{\,n}\,\right\| &\,\leq\, \sum\limits_{v \,=\, q}^{\,p \,-\, 1}\,\left\|\,x_{\,v} \,-\, x_{\,v \,+\, 1},\, b_{\,2},\, \cdots,\, b_{\,n}\,\right\|\\
&=\, \sum\limits_{v \,=\, p}^{\,p \,-\, 1}\,\left\|\,T^{\,v}\,x_{\,0} \,-\, T^{\,v}\,x_{\,1},\, b_{\,2},\, \cdots,\, b_{\,n}\,\right\|\\
& \,\leq\, \sum\limits_{v \,=\, q}^{\,p \,-\, 1}\,a_{\,v}\,\left\|\,x_{\,0} \,-\, x_{\,1},\, b_{\,2},\, \cdots,\, b_{\,n}\,\right\|\; \;[\;\text{by}\;(\,\ref{eq1.8}\,)\;]
\end{align*}
If \,$x_{\,0} \,=\, x_{\,1}$, then a fixed point is obtained.\,Let \,$x_{\,0} \,\neq\, x_{\,1}$\, and \,$r$\, be a positive integer with \,$r \,>\, \left\|\,x_{\,0} \,-\, x_{\,1},\, b_{\,2},\, \cdots,\, b_{\,n}\,\right\|$.\,As the series \,$\sum\limits_{k \,=\, 1}^{\,\infty}\,a_{\,k}$\, is convergent, for \,$\epsilon \,>\, 0$\, arbitrary there exists a positive integer \,$k_{\,0}$\, such that \,$\sum\limits_{v \,=\, q}^{\,p \,-\, 1}\,a_{\,v} \,<\, \dfrac{\epsilon}{r}$\, if \,$p \,>\, q \,\geq\, k_{\,0}$.\,Then for    \,$p \,>\, q \,\geq\, k_{\,0}$, we have
\[\left\|\,x_{\,p} \,-\, x_{\,q},\, b_{\,2},\, \cdots,\, b_{\,n}\,\right\| \,<\, \dfrac{\epsilon}{r}\,\left\|\,x_{\,0} \,-\, x_{\,1},\, b_{\,2},\, \cdots,\, b_{\,n}\,\right\| \,<\, \epsilon.\]
Therefore \,$\left\{\,x_{\,k}\,\right\}$\, is a \,$b$-Cauchy sequence and \,$b$-completeness of \,$X$\, implies the existence of \,$\xi \,\in\, X$\, such that \,$\lim\limits_{k \,\to\, \infty}\,x_{\,k} \,=\, \xi$.\,If \,$k$\, be a positive integer then
\begin{align*}
&\left\|\,\xi \,-\, T\,\xi,\, b_{\,2},\,\cdots,\,b_{\,n}\,\right\| \,\leq\, \left\|\,\xi \,-\, x_{\,k \,+\, 1},\, b_{\,2},\,\cdots,\,b_{\,n}\,\right\| \,+\, \left\|\,x_{k \,+\, 1} \,-\, T\,\xi,\, b_{\,2},\,\cdots,\,b_{\,n}\,\right\|\\
&=\, \left\|\,\xi \,-\, x_{\,k \,+\, 1},\, b_{\,2},\,\cdots,\,b_{\,n}\,\right\| \,+\, \left\|\,T\,x_{\,k} \,-\, T\,\xi,\, b_{\,2},\,\cdots,\,b_{\,n}\,\right\|\\
&\leq\, \left\|\,\xi \,-\, x_{\,k \,+\, 1},\, b_{\,2},\,\cdots,\,b_{\,n}\,\right\| \,+\, a_{\,1}\left\|\,x_{\,k} \,-\, \xi,\, b_{\,2},\,\cdots,\,b_{\,n}\,\right\|\; \;[\;\text{by (\ref{eq1.8})}\;]\\
&\,\to\, 0\, \;\text{as}\, \,k \,\to\, \infty
\end{align*} 
So, \,$T\,\xi \,=\, \xi$\, and \,$\xi$\, becomes a fixed point of \,$T$.\,We now prove the uniqueness of \,$\xi$.\,If \,$\eta$\, be a fixed point of \,$T$\, then for any positive integer \,$k$, \,$\eta \,=\, T^{\,k}\,\eta$\, and \,$\xi \,=\, T^{\,k}\,\xi$.\,So,
\[\left\|\,\xi \,-\, \eta,\, b_{\,2},\, \cdots,\, b_{\,n}\,\right\| \,=\, \left\|\,T^{\,k}\,\xi \,-\, T^{\,k}\,\eta,\, b_{\,2},\, \cdots,\, b_{\,n}\,\right\| \,\leq\, a_{\,k}\,\left\|\,\xi \,-\, \eta,\, b_{\,2},\, \cdots,\, b_{\,n}\,\right\|.\]
If  \,$\left\|\,\xi \,-\, \eta,\, b_{\,2},\, \cdots,\, b_{\,n}\,\right\| \,>\, 0$\, then \,$a_{\,k \,} \,\geq\, 1$\, for all \,$k$.\,So, \,$a_{\,k}$\, cannot tends to zero and this contradiction shows that \,$\xi \,=\, \eta$\, and the theorem is proved.
\end{proof}

\begin{remark}
We now deduce Theorem \ref{th1} from the above Theorem \ref{thmm1}.
\end{remark}

\begin{proof}
Since \,$T$\, is a \,$b$-contraction mapping, there exists \,$0 \,<\, \alpha \,<\, 1$\, such that
\begin{equation}\label{eq1.9}
\left\|\,T\,x \,-\, T\,y,\, b_{\,2},\,\cdots,\,b_{\,n}\,\right\| \,\leq\, \alpha\,\left\|\,x \,-\, y,\, b_{\,2},\, \cdots,\, b_{\,n}\,\right\|\; \;\forall\; x,\,y \,\in\, X.
\end{equation}
For \,$x,\,y \,\in\, X$, we get from (\ref{eq1.9})
\begin{align*}
\left\|\,T^{\,2}\,x \,-\, T^{\,2}\,y,\, b_{\,2},\,\cdots,\,b_{\,n}\,\right\| &\,\leq\, \alpha\,\left\|\,T\,x \,-\, T\,y,\, b_{\,2},\, \cdots,\, b_{\,n}\,\right\|\\
& \,\leq\, \alpha^{\,2}\,\left\|\,x \,-\, y,\, b_{\,2},\, \cdots,\, b_{\,n}\,\right\|.
\end{align*}
\begin{align*}
\left\|\,T^{\,3}\,x \,-\, T^{\,3}\,y,\, b_{\,2},\,\cdots,\,b_{\,n}\,\right\| &\,\leq\, \alpha\,\left\|\,T^{\,2}\,x \,-\, T^{\,2}\,y,\, b_{\,2},\, \cdots,\, b_{\,n}\,\right\|\\
& \,\leq\, \alpha^{\,3}\,\left\|\,x \,-\, y,\, b_{\,2},\, \cdots,\, b_{\,n}\,\right\|
\end{align*}
and in general 
\[\left\|\,T^{\,k}\,x \,-\, T^{\,k}\,y,\, b_{\,2},\,\cdots,\,b_{\,n}\,\right\| \,\leq\, \alpha^{\,k}\,\left\|\,x \,-\, y,\, b_{\,2},\, \cdots,\, b_{\,n}\,\right\|.\]
Since the series \,$\sum\limits_{k \,=\, 1}^{\,\infty}\,\alpha^{\,k}$\, is convergent, by previous Theorem \,$T$\, has a fixed point in \,$X$.
\end{proof}

\begin{theorem}\label{th1.2}
Let \,$T \,:\, X \,\to\, X$\, be an operator, where \,$X$\, is a \,$b$-complete linear \,$n$-normed space and \,$T$\, satisfies the condition 
\begin{align}
&\left\|\,T\,x \,-\, T\,y,\, b_{\,2},\,\cdots,\,b_{\,n}\,\right\|\nonumber\\
& \,\leq\, \beta\,\left\{\,\left\|\,x \,-\, T\,x,\, b_{\,2},\, \cdots,\, b_{\,n}\,\right\| \,+\, \left\|\,y \,-\, T\,y,\, b_{\,2},\, \cdots,\, b_{\,n}\,\right\|\,\right\}\label{eq1.21}
\end{align}
where \,$0 \,<\, \beta \,<\, \dfrac{1}{2}$\, and \,$x,\,y \,\in\, X$.\,Then there exists a unique fixed point \,$x_{\,0}$\, of the operator \,$T$\, in \,$X$\, provided the set \,$\left\{\,x_{\,0},\, b_{\,2},\, \cdots,\, b_{\,n}\,\right\}$\, is linearly independent.
\end{theorem}

\begin{proof}
Let \,$x \,\in\, X$\, and \,$x_{\,1} \,=\, T\,x,\, x_{\,2} \,=\, T\,x_{\,1},\, x_{\,3} \,=\, T\,x_{\,2},\, \cdots,\, ,\, x_{\,k} \,=\, T\,x_{\,k \,-\, 1},\, \cdots$. Then
\begin{align*}
&\left\|\,x_{\,1} \,-\, x_{\,2},\, b_{\,2},\, \cdots,\, b_{\,n}\,\right\| \,=\, \left\|\,T\,x \,-\, T\,x_{\,1},\, b_{\,2},\,\cdots,\,b_{\,n}\,\right\|\\
&\leq\, \beta\left\{\,\left\|\,x \,-\, T\,x,\, b_{\,2},\, \cdots,\, b_{\,n}\,\right\| \,+\, \left\|\,x_{\,1} \,-\, T\,x_{\,1},\, b_{\,2},\, \cdots,\, b_{\,n}\,\right\|\,\right\}\\
&=\, \beta\left\{\,\left\|\,x \,-\, T\,x,\, b_{\,2},\, \cdots,\, b_{\,n}\,\right\| \,+\, \left\|\,x_{\,1} \,-\, x_{\,2},\, b_{\,2},\, \cdots,\, b_{\,n}\,\right\|\,\right\}.\\
&\Rightarrow\left\|\,x_{\,1} \,-\, x_{\,2},\, b_{\,2},\, \cdots,\, b_{\,n}\,\right\| \,\leq\, \dfrac{\beta}{1 \,-\, \beta}\,\left\|\,x \,-\, T\,x,\, b_{\,2},\, \cdots,\, b_{\,n}\,\right\|.
\end{align*}
And
\begin{align*}
&\left\|\,x_{\,2} \,-\, x_{\,3},\, b_{\,2},\, \cdots,\, b_{\,n}\,\right\| \,=\, \left\|\,T\,x_{\,1} \,-\, T\,x_{\,2},\, b_{\,2},\,\cdots,\,b_{\,n}\,\right\|\\
&\leq\, \beta\left\{\,\left\|\,x_{\,1} \,-\, T\,x_{\,1},\, b_{\,2},\, \cdots,\, b_{\,n}\,\right\| \,+\, \left\|\,x_{\,2} \,-\, T\,x_{\,2},\, b_{\,2},\, \cdots,\, b_{\,n}\,\right\|\,\right\}\\
&=\, \beta\left\{\,\left\|\,x_{\,1} \,-\, x_{\,2},\, b_{\,2},\, \cdots,\, b_{\,n}\,\right\| \,+\, \left\|\,x_{\,2} \,-\, x_{\,3},\, b_{\,2},\, \cdots,\, b_{\,n}\,\right\|\,\right\}.\\
&\Rightarrow\left\|\,x_{\,2} \,-\, x_{\,3},\, b_{\,2},\, \cdots,\, b_{\,n}\,\right\| \,\leq\, \dfrac{\beta}{1 \,-\, \beta}\,\left\|\,x_{\,1} \,-\, x_{\,2},\, b_{\,2},\, \cdots,\, b_{\,n}\,\right\|\\\
&\hspace{4.7cm}\,\leq\, \left(\,\dfrac{\beta}{1 \,-\, \beta}\,\right)^{\,2}\,\left\|\,x \,-\, T\,x,\, b_{\,2},\, \cdots,\, b_{\,n}\,\right\|.
\end{align*}
Similarly,
\[\left\|\,x_{\,3} \,-\, x_{\,4},\, b_{\,2},\, \cdots,\, b_{\,n}\,\right\| \,\leq\, \left(\,\dfrac{\beta}{1 \,-\, \beta}\,\right)^{\,3}\,\left\|\,x \,-\, T\,x,\, b_{\,2},\, \cdots,\, b_{\,n}\,\right\|.\]
In general, if \,$k$\, is any positive integer, then
\[\left\|\,x_{\,k} \,-\, x_{\,k\,+\,1},\, b_{\,2},\, \cdots,\, b_{\,n}\,\right\| \,\leq\, \left(\,\dfrac{\beta}{1 \,-\, \beta}\,\right)^{\,k}\,\left\|\,x \,-\, T\,x,\, b_{\,2},\, \cdots,\, b_{\,n}\,\right\|.\]
So, if \,$p$\, is any positive integer, then 
\begin{align}
&\left\|\,x_{\,k} \,-\, x_{\,k \,+\, p},\, b_{\,2},\, \cdots,\, b_{\,n}\,\right\|\nonumber\\
&\leq\,\left\|\,x_{\,k} - x_{\,k \,+\, 1},\, b_{\,2},\, \cdots,\, b_{\,n}\,\right\| \,+\, \left\|\,x_{\,k \,+\, 1} - x_{\,k \,+\, 2},\, b_{\,2},\, \cdots,\, b_{\,n}\,\right\|\,+\,\cdots\nonumber\\
&\hspace{1cm}\cdots\,+\,\left\|\,x_{\,k\,+\,p\,-1} - x_{\,k \,+\, p},\, b_{\,2},\, \cdots,\, b_{\,n}\,\right\|\nonumber\\
&\leq\, \left(\,r^{\,k} \,+\, r^{\,k\,+\,1} \,+\, \cdots\,+\, r^{\,k\,+\,p\,-\,1}\,\right)\,\left\|\,x \,-\, T\,x,\, b_{\,2},\, \cdots,\, b_{\,n}\,\right\|, \;\text{where}\;\,r \,=\, \dfrac{\beta}{1 \,-\, \beta}\nonumber\\
&<\, \dfrac{r^{\,k}}{1 \,-\, r}\,\left\|\,x \,-\, T\,x,\, b_{\,2},\, \cdots,\, b_{\,n}\,\right\|.\label{eq1.3}
\end{align}
Since \,$0 \,<\, \beta \,<\, \dfrac{1}{2}$, we have \,$0 \,<\, r \,<\, 1$\, and so by (\ref{eq1.3}), \,$\left\|\,x_{\,k} \,-\, x_{\,k \,+\, p},\, b_{\,2},\, \cdots,\, b_{\,n}\,\right\| \,\to\, 0$\, as \,$k \,\to\, \infty$.\,Therefore, \,$\left\{\,x_{\,k}\,\right\}$\, is a \,$b$-Cauchy sequence.\,Since \,$X$\, is \,$b$-complete linear \,$n$-normed space, the sequence \,$\left\{\,x_{\,k}\,\right\}$\, is convergent in the semi-normed space \,$\left(\,X,\, \left\|\,\cdot,\,b_{\,2},\, \cdots,\, b_{\,n}\,\right\|\,\right)$.\;So, let \,$\lim\limits_{k \,\to\, \infty}\,x_{\,k} \,=\, x_{\,0}$\, with the property that \,$\left\{\,x_{\,0},\, b_{\,2},\, \cdots,\, b_{\,n}\,\right\}$\, is linearly independent.\,We now show that \,$x_{\,0}$\, is a fixed point of \,$T$.\,We have
\begin{align*}
&\left\|\,x_{\,0} \,-\, T\,x_{\,0},\, b_{\,2},\, \cdots,\, b_{\,n}\,\right\| \,\leq\, \left\|\,x_{\,0} \,-\, x_{\,k},\, b_{\,2},\, \cdots,\, b_{\,n}\,\right\| \,+\, \left\|\,x_{\,k} \,-\, T\,x_{\,0},\, b_{\,2},\, \cdots,\, b_{\,n}\,\right\|\\
&=\, \left\|\,x_{\,0} \,-\, x_{\,k},\, b_{\,2},\, \cdots,\, b_{\,n}\,\right\| \,+\, \left\|\,T\,x_{\,k\,-\,1} \,-\, T\,x_{\,0},\, b_{\,2},\, \cdots,\, b_{\,n}\,\right\|\\
&\leq\,  \beta\left\{\,\left\|\,x_{k\,-\,1} - T\,x_{k\,-\,1},\, b_{\,2},\, \cdots,\, b_{\,n}\,\right\| + \left\|\,x_{\,0} - T\,x_{\,0},\, b_{\,2},\, \cdots,\, b_{\,n}\,\right\|\,\right\}\,+\\
&\hspace{1cm}+\,\left\|\,x_{\,0} - x_{\,k},\, b_{\,2},\, \cdots,\, b_{\,n}\,\right\|.
\end{align*}
This implies that
\begin{align*}
&(\,1 \,-\, \beta\,)\,\left\|\,x_{\,0} - T\,x_{\,0},\, b_{\,2},\, \cdots,\, b_{\,n}\,\right\|\\
& \,\leq\, \left\|\,x_{\,0} - x_{\,k},\, b_{\,2},\, \cdots,\, b_{\,n}\,\right\| \,+\, \beta\,\left\|\,x_{k\,-\,1} - x_{\,k},\, b_{\,2},\, \cdots,\, b_{\,n}\,\right\|\\    
&\to\, 0\; \;\text{as}\; \,k \,\to\, \infty\; \;\left[\;\text{since}\;\lim\limits_{k \,\to\, \infty}\,x_{\,k} \,=\, x_{\,0}\,\right].
\end{align*}
So, \,$T\,x_{\,0} \,=\, x_{\,0}$\, and \,$x_{\,0}$\, is a fixed point of \,$T$.\,We now prove that \,$x_{\,0}$\, is the only fixed point of \,$T$.\,Let \,$T\,y_{\,0} \,=\, y_{\,0}$\, such that \,$\left\{\,x_{\,0} \,-\, y_{\,0},\, b_{\,2},\, \cdots,\, b_{\,n}\,\right\}$\, is linearly independent.\,Then
\begin{align*}
&\left\|\,x_{\,0} \,-\, y_{\,0},\, b_{\,2},\, \cdots,\, b_{\,n}\,\right\| \,=\, \left\|\,T\,x_{\,0} \,-\, T\,y_{\,0},\, b_{\,2},\, \cdots,\, b_{\,n}\,\right\| \\
&\leq\, \beta\left\{\,\left\|\,x_{\,0} \,-\, T\,x_{\,0},\, b_{\,2},\, \cdots,\, b_{\,n}\,\right\| \,+\, \left\|\,y_{\,0} \,-\, T\,y_{\,0},\, b_{\,2},\, \cdots,\, b_{\,n}\,\right\|\,\right\} \,=\, 0.
\end{align*}
This shows that \,$x_{\,0} \,=\, y_{\,0}$.\,This proves the theorem.
\end{proof}

\end{document}